\documentclass[10pt]{amsart}
\usepackage{mathrsfs}
\usepackage{amsfonts} 
\usepackage{graphicx,latexsym,bm,amsmath,amssymb,verbatim,multicol,lscape}
\newtheorem{thm}{Theorem} [section]
\newtheorem{cor}[thm]{Corollary}

\theoremstyle{definition}

\theoremstyle{remark}

\numberwithin{equation}{section}

\begin{document}
\title [Asymptotic behavior of LCM of consecutive arithmetic progression terms]
{Asymptotic behavior of the least common multiple of consecutive arithmetic progression terms}

\author{Guoyou Qian}
\address{Center for Combinatorics, Nankai University, Tianjin 300071, P.R. China}
\email{qiangy1230@163.com,qiangy1230@gmail.com}
\author{Shaofang Hong}
\address{Yangtze Center of Mathematics,
Sichuan University, Chengdu 610064, P.R. China and Mathematical
College, Sichuan University, Chengdu 610064, P.R. China}
\email{sfhong@scu.edu.cn, s-f.hong@tom.com, hongsf02@yahoo.com }

\thanks{Hong is the corresponding author and was supported partially by
National Science Foundation of China Grant \# 10971145 and by the Ph.D.
Programs Foundation of Ministry of Education of China Grant \#20100181110073}
\date{\today}%
\begin{abstract}
Let $l$ and $m$ be two integers with $l>m\ge 0$, and let $a$ and $b$
be integers with $a\ge 1$ and $a+b\ge 1$.  In this paper, we prove
that $\log {\rm lcm}_{mn<i\le ln}\{ai+b\} =An+o(n)$, where $A$ is
a constant depending on $l, m$ and $a$.
\end{abstract}

\maketitle

\section{Introduction}
Chebyshev \cite{[Ch]} initiated the study  of the  least common
multiple of consecutive positive integers  for the first significant
attempt to prove prime number theorem. Actually, Chebyshev
introduced the function $\Psi(x)=\sum_{p^k\le x}\log p= \log {\rm
lcm}_{1\le i\le x}\{i \}$, and then obtained an equivalent form of
the prime number theorem, which states that $\log {\rm lcm}(1, ...,
n)\sim n$ as $n$ goes to infinity. Since then, this topic received
many authors' attention. Hanson \cite{[Ha]} and Nair \cite{[N]} got
upper and lower bounds of ${\rm lcm}_{1\le i\le n}\{i\}$,
respectively. Farhi \cite{[F]} investigated the least common
multiple of arithmetic progressions, while Farhi and Kane
\cite{[FK]} and Hong and Yang \cite{[HY]} studied the least common
multiple of consecutive positive integers. Hong and Qian \cite{[HQ]}
obtained some results on the least common multiple of consecutive
arithmetic progression terms. Let $a$ and $b$ be integers such that
$a\ge 1$, $a+b\ge 1$ and $\gcd(a, b)=1$. On the other hand, Bateman,
Kalb and Stenger \cite{[BKS]} proved that
$$
\log{\rm lcm}_{1\le i\le n}\{ai+b\}\sim
\frac{an}{\varphi(a)}\sum_{r=1 \atop \gcd(r, a)=1}^a\frac{1}{r}
$$
as $n\rightarrow\infty $, where $\varphi(a)$ denotes the number of
integers relatively prime to $a$ between 1 and $a$. Farhi \cite{[F]}
showed that ${\rm lcm}_{1\le i\le n}\{i^2+1\}\ge 0.32\cdot
(1.442)^n$ for all $n\ge 1$. Qian, Tan and Hong \cite{[QTH]} proved
that for any given positive integer $k$, $\log{\rm lcm}_{0\le i\le
k}\{(n+i)^2+1\} \sim 2(k+1)\log n$ as $n\rightarrow \infty$. 

In this paper, we mainly focus on the least common multiple $${\rm
lcm}_{mn<i\le ln}\{ai+b\}$$ of consecutive arithmetic progression
terms, where $l$ and $m$ are integers such that $l>m\ge 0$ and
$a\ge 1$ and $b$ are integers such that $a+b\ge 1$. Evidently, if
$\gcd(a, b)=d$, then
\begin{align*} 
\log {\rm lcm}_{mn<i\le ln}\{ai+b\}&=\log {\rm lcm}_{mn<i\le
ln}\{a_1i+b_1\}+\log d\\
&=\log {\rm lcm}_{mn<i\le ln}\{a_1i+b_1\}+O(1),
\end{align*}
where $a_1=a/d$ and $b_1=b/d$.  So it is sufficient to consider the
case $\gcd(a, b)=1$. Let $h$ and $k$ be any given two relatively
prime positive integers. We define

\begin{align}\label{eq:1.1}
\vartheta(x;h,k):=\sum_{{\rm prime}\ p\le x \atop p\equiv k\pmod
{h}}\log p.
\end{align}
Then the prime number theorem for arithmetic progressions states
that
\begin{align}\label{eq:1.2}\vartheta(x;h,k)=\frac{x}{\varphi(h)}
+O\big(x\exp(-c_1\sqrt{\log x})\big),
\end{align}
where $c_1>0$ is a constant (see, for example, \cite{[MV]}).

As usual, for any prime  $p$, we let $v_{p}$ be the normalized
$p$-adic valuation of $\mathbb{N}^*$, i.e., $v_p(a)=s$ if
$p^{s}\parallel a$.  For any real number $x$, we let $\lfloor
x\rfloor$ denote the largest integer no more than $x$.
We can now state the main result of this paper.\\

\begin{thm}\label{mainresult}  Let $l$ and $m$ be integers with $l>m\ge 0$
and let $a\ge 1$ and $b$ be integers such that $a+b\ge 1$ and
$\gcd(a, b)=1$. Then
$$\log {\rm lcm}_{mn<i\le ln}\{ ai+b\}
\sim \frac{an}{\varphi(a)}\sum_{r=1\atop \gcd(r,a)=1}^{a}A_r
$$
as $n\rightarrow \infty$, where
\begin{align}\label{eq:1.3} A_r:={\left\{
\begin{array}{rl}
\frac{l}{r}, &\text{if} \ l\ge \frac{(a+r)m}{r},\\
\sum_{i=0}^{K_r-1}\frac{l-m}{r+ai}+\frac{l}{r+aK_r}, &\text{if} \
m<l<\frac{(a+r)m}{r}
\end{array}
\right.}
\end{align}
with
\begin{align}\label{eq:1.4}
K_r:=\Big\lfloor\frac{al-(l-m)r}{a(l-m)}\Big\rfloor.
\end{align}
\end{thm}

From Theorem \ref{mainresult}, we can deduce immediately the following interesting results.\\

\begin{cor}\label{cor1}  Let $l$ and $m$ be integers with $l>m\ge 0$ and let $a$
and $b$ be integers such that $a\ge 1$, $a+b\ge 1$ and $\gcd(a,
b)=1$. If $l\ge (a+1)m$,  then
$$\log {\rm lcm}_{mn<i\le ln}\{ ai+b\}
\sim  \frac{aln}{\varphi(a)}\sum_{r=1\atop \gcd(r,a)=1}^{a}
\frac{1}{r}
$$
as $n\rightarrow \infty$.
\end{cor}

Taking $l=1$ and  $m=0$, Corollary \ref{cor1} then becomes the
Bateman-Kalb-Stenger theorem \cite{[BKS]}. \\

\begin{cor}  Let $l$ and $m$ be integers with $l>m\ge 0$. Then
$$\log {\rm lcm}_{mn<i\le ln}\{i\}\sim Bn
$$
as $n\rightarrow \infty$, where
$$
B={\left\{
\begin{array}{rl}
l, &{\it if} \ l\ge 2m,\\
\displaystyle\frac{l}{\lfloor\frac{l}{l-m}\rfloor}
+(l-m)\sum_{i=1}^{\lfloor\frac{m}{l-m}\rfloor} \frac{1}{i}, & {\it
if} \ m<l<2m.
\end{array}
\right.}
$$
\end{cor}

The next section will devote to the proof of Theorem
\ref{mainresult}.

\section{Proof of Theorem \ref{mainresult}}

In this section, we show Theorem \ref{mainresult}.\\
\\
{\it Proof of Theorem \ref{mainresult}.} For simplicity,  we let
$$L_{m, l}(n):={\rm lcm}_{mn<i\le ln}\{b+ai\}$$
and by  $P_{m, l}(n)$ we denote the set of all the prime factors of
$L_{m, l}(n)$. Denote by $R(a)$ the set of all the integers
relatively prime to $a$ between $1$ and $a$.

For any given integer $b$ with $a+b\ge 1$ and $\gcd(a, b)=1$, we can
assume that $b=b_0+qa$ for some integers $q$ and $b_0\in R(a)$. Then
\begin{align*}
\log{\rm lcm}\big(b+a(mn+1),..., b+aln\big)&=\log{\rm lcm}\big(b_0+a(mn+1+q),..., b_0+a(ln+q)\big)\\
&=\log{\rm lcm}_{mn<i\le ln}\{b_0+ai\}+O\big(\log n\big)\\
&=\log{\rm lcm}_{mn<i\le ln}\{b_0+ai\}+o(n).
\end{align*}
So to prove Theorem \ref{mainresult},  we may assume that $b\in
R(a)$ in what follows.

First of all, we have
\begin{align*}
L_{m, l}(n)=\prod_{p\in P_{m, l}(n)} p^{v_p(L_{m,
l}(n))}=\Big(\prod_{p\in P_{m, l}(n)} p\Big)\Big(\prod_{p\in P_{m,
l}(n)\atop p^2|L_{m, l}(n)} p^{v_p(L_{m, l}(n))-1}\Big).
\end{align*}
So
\begin{equation}\label{eq:2.1}
\log L_{m, l}(n)=\sum_{p\in P_{m, l}(n)}\log p+\sum_{p\in P_{m,
l}(n)\atop p^2|L_{m, l}(n)}\big(v_p(L_{m, l}(n))-1\big)\log p.
\end{equation}
If $p^2|L_{m, l}(n)$, then $p^2|(b+ai_0)$ for some integer
$mn<i_0\le ln$, which implies that $p^2\le b+ai_0\le b+aln$.  Hence
for any  prime  $p$ running over the second summation in
(\ref{eq:2.1}), we have $p\le \sqrt {b+aln}$. Since  $p^{v_p(L_{m,
l}(n))}\le b+aln$ for any $p\in P_{m, l}(n)$, one has
$$v_p(L_{m, l}(n))\le \frac{\log (b+aln)}{\log p}.$$
It follows from the prime number theorem that
\begin{align*}
\sum_{p\in P_{m, l}(n)\atop p^2|L_{m, l}(n)} \big(v_p(L_{m,
l}(n))-1\big) \log p &\le \sum_{p\le \sqrt{b+aln}}\frac{\log
(b+aln)}{\log p}\log p= \sum_{p\le
\sqrt{b+aln}}\log (b+aln)\\
&\ll \frac{\sqrt{b+aln}}{\log \sqrt {b+aln}}\log(b+aln)=2\sqrt{
b+aln}.
\end{align*}
Then by (\ref{eq:2.1}) we obtain
\begin{equation}\label{eq:2.2}
\log L_{m, l}(n)=\sum_{p\in P_{m, l}(n)}\log p+O(\sqrt{n}).
\end{equation}

To estimate $\sum_{p\in P_{m, l}(n)}\log p$, we first need to find a
characterization on the primes in the set $P_{m, l}(n)$. Since
$\gcd(a, b)=1$, we get that $\gcd(a, b+ai)=1$ for all positive
integers $i$, which implies that each prime $p$ in $P_{m, l}(n)$ is
relatively prime to $a$.  Then each prime $p\in P_{m, l}(n)$ is
congruent to some $r\in R(a)$ modulo $a$. It then follows from
(\ref{eq:2.2}) that
\begin{equation}\label{eq:2.3}
\log L_{m, l}(n)=\sum_{r\in R(a)}\sum_{p\in P_{m,l}(n)\atop p\equiv
r\pmod a}\log p+O(\sqrt{n}).
\end{equation}
By (\ref{eq:2.3}), to estimate $\log L_{m, l}(n)$, it suffices to
estimate the sum
\begin{equation}\label{eq:2.4}
\sum_{p\in P_{m,l}(n)\atop p\equiv r\pmod a}\log p
\end{equation}
for each $r\in R(a)$, which will be done in the following.

Now fix an $r\in R(a)$. Then there exists only one $r'\in R(a)$ such
that $rr'\equiv b\pmod a$. If $p$ is a prime with $p\equiv r\pmod
a$, then $r'p\equiv b\pmod a$. So any term divisible by $p$ in the
arithmetic progression $\{b+ai\}_{i=1}^{\infty}$ must be of the form
$(r'+aj)p$ with $j\ge 0$ being an integer. Thus for any prime
$p\equiv r\pmod a$,  we have that $p\in P_{m, l}(n)$  if and only if
there exists a nonnegative integer $i_0$ such that
$$b+amn<(r'+ai_0)p\le b+aln.$$
In other words, a prime $p$ congruent to $r$ modulo $a$ is in $P_{m,
l}(n)$ if and only if
$$\frac{b+amn}{r'+ai_0}<p\le \frac{b+aln}{r'+ai_0}$$
for some nonnegative integer $i_0$. Therefore
\begin{align}\label{eq:2.5}
\Big\{p\in P_{m, l}(n): p\equiv r\ ({\rm mod}\
a)\Big\}=\bigcup_{i=0}^{\infty} \Big\{{\rm prime}\ q\equiv r\ ({\rm
mod}\ a): \frac{b+amn}{r'+ai}<q\le \frac{b+aln}{r'+ai}\Big\}.
\end{align}
In what follows we transfer the union in the right-hand side of
(\ref{eq:2.5}) into a union of finitely many sets.

If $p\not|a$, then exactly one term of any consecutive $p$ terms in
the arithmetic progression $\{b+ai\}_{i=1}^{\infty}$ is divisible by
$p$. Therefore, for any prime $p$ with $p\le (l-m)n$ and $p\not|a$,
there is at least one integer $i_0$ such that $mn<i_0\le ln$ and
$b+ai_0\equiv 0\pmod p$. Hence the primes $p$ satisfying $p\le
(l-m)n$ and $p\not|a$ are all in the set $P_{m, l}(n)$.

For convenience, we let $H=K_{r'}$, where $K_{r'}$ is defined as in
(\ref{eq:1.4}). Then we have
\begin{equation}\label{eq:2.6}
H=\Big\lfloor \frac{al-(l-m)r'}{a(l-m)}\Big\rfloor.
\end{equation}
Thus for any positive integer $n$, one has
\begin{equation}\label{eq:2.7}
(l-m)n<\frac{b+aln}{r'+aH}.
\end{equation}and for sufficiently large $n$, we have
$$\frac{b+aln}{r'+a(H+1)}<(l-m)n.
$$
In the rest of the proof, we always assume that $n$ is a sufficient
large integer.

Since $\frac{b+aln}{r'+ai}$ strictly decreases as $i$ increases, it
follows immediately that for any integer $i>H$, we have
$$(l-m)n>\frac{b+aln}{r'+ai}>\frac{b+amn}{r'+ai},$$
which implies that
$$\Big(\frac{b+amn}{r'+ai}, \frac{b+aln}{r'+ai}
\Big]\subseteq \big(0, (l-m)n\big]$$ for any integer $i>H$.
Therefore
\begin{align*}
&\bigcup_{i=H+1}^{\infty}\Big\{{\rm prime}\ q\equiv r \ ({\rm mod}\
a): \frac{b+amn}{r'+ai}<q\le \frac{b+aln}{r'+ai}\Big\}\\
&\quad\subseteq \Big\{{\rm prime}\ q\equiv r \ ({\rm mod}\ a): q\le
(l-m)n\Big\}.
\end{align*}
Hence by (\ref{eq:2.5}),
\begin{align}\label{eq:2.8}
\big\{p\in P_{m, l}(n): p\equiv r\ ({\rm mod}\ a)\big\}=\
&\bigcup_{i=0}^{H}\Big\{{\rm prime}\ q\equiv r \ ({\rm mod}\ a):
\frac{b+amn}{r'+ai}<q\le \frac{b+aln}{r'+ai}\Big\}\\
\nonumber&\quad \ \bigcup \Big\{{\rm prime}\ q\equiv r \ ({\rm mod}\
a): q\le (l-m)n\Big\}.
\end{align}

To estimate the sum (\ref{eq:2.4}), we consider the following two cases.\\

{\bf Case 1.}  $l>\frac{(a+r')m}{r'}$. Then $(a+r')m< lr'$. It
follows from (\ref{eq:2.6}) that $H=0$. So by (\ref{eq:2.8}), we
only need to deal with the following union
\begin{align}\label{eq:2.9}
\Big\{{\rm prime}\ q\equiv r \ ({\rm mod}\ a):&
\frac{b+amn}{r'}<q\le \frac{b+aln}{r'}\Big\}\bigcup \Big\{{\rm
prime}\ q\equiv r \ ({\rm mod}\ a): q\le (l-m)n\Big\},
\end{align}
which is equal to the set $\big\{p\in P_{m, l}(n): p\equiv r\ ({\rm
mod}\ a)\big\}$.

Since $(a+r')m< lr'$, we have
$$b+amn< (l-m)r'n+b\le b+aln,$$
which implies that
$$\frac{b+amn}{r'}< (l-m)n+\frac{b}{r'}\le \frac{b+aln}{r'}.$$
Then the union of the following three intervals
$$\big(0, (l-m)n\big] \bigcup \big((l-m)n, (l-m)n+\frac{b}{r'} \big]
\bigcup\big(\frac{b+amn}{r'}, \frac{b+aln}{r'}\big]$$ equals the
interval $(0, \frac{b+aln}{r'}]$. Hence by (\ref{eq:2.9}), we obtain
that all the primes $p\le \frac{b+aln}{r'}$ congruent to $r$ modulo
$a$ and not belonging to the interval $\big((l-m)n,
(l-m)n+\frac{b}{r'} \big]$ are contained in the set $\big\{p\in
P_{m, l}(n): p\equiv r\ ({\rm mod}\ a)\big\}$. So from
(\ref{eq:1.1}) and (\ref{eq:1.2}), we derive that if
$l>\frac{(a+r')m}{r'}$, then
\begin{align}\label{eq:2.10}
\sum_{p\in P_{m,l}(n)\atop p\equiv r\pmod a}\log p &=\sum_{{\rm
prime}\ p\equiv r({\rm mod}\ a)\atop p \le \frac{b+aln}{r'}}\log p
+O\Big(\sum_{{\rm prime}\ p\equiv r({\rm mod}\ a)\atop
(l-m)n<p\le (l-m)n+\frac{b}{r'}}\log p\Big)\\
\nonumber&=\vartheta\Big( \frac{b+aln}{r'}; a, r\Big)+ O\Big(\log
n\Big)\\
\nonumber&=
\frac{a}{\varphi(a)}\frac{ln}{r'}+O\Big(n\exp(-c_1\sqrt{\log
n})\Big).
\end{align}

{\bf Case 2.} $l\le \frac{(a+r')m}{r'}$.  In this case, $H\ge 1$.
Clearly, $H>\frac{al-(l-m)r'}{a(l-m)}-1$ by (\ref{eq:2.6}). So we
obtain
\begin{align*}
(l-m)n-\frac{b+amn}{r'+aH}&=\frac{(l-m)r'n+aH(l-m)n-(b+amn)}{r'+aH}\\
&>
\frac{(l-m)r'n+a(l-m)(\frac{al-(l-m)r'}{a(l-m)}-1)n-(b+amn)}{r'+aH}\\
&=\frac{-b}{r'+aH}.
\end{align*}
It follows from $b\in R(a)$ and $r'\in R(a)$ that
\begin{equation}\label{eq:2.11}
\frac{b+amn}{r'+aH}< (l-m)n+\frac{b}{r'+aH}< (l-m)n+1.
\end{equation}
On the other hand, it is obvious that
\begin{align*}
&\Big\{{\rm prime}\ q\equiv r \ ({\rm mod}\ a): q\le
(l-m)n\Big\}=\Big\{{\rm prime}\ q\equiv r \ ({\rm mod}\ a):
q<(l-m)n+1\Big\}.
\end{align*}
 Thus we derive from (\ref{eq:2.7}) and
(\ref{eq:2.11}) that
\begin{align*}
&\Big\{{\rm prime}\ q\equiv r \ ({\rm mod}\ a): q\le
(l-m)n\Big\}\bigcup\Big\{{\rm
prime}\ q\equiv r \ ({\rm mod}\ a): \frac{b+amn}{r'+aH}<q\le \frac{b+aln}{r'+aH}\Big\}\\
&=\Big\{{\rm prime}\ q\equiv r({\rm mod}\ a): q\le
\frac{b+aln}{r'+aH}\Big\}.
\end{align*}
Then by (\ref{eq:2.8}), we have
\begin{equation}\label{eq:2.12}
\big\{{\rm prime}\ p\in P_{m, l}(n): p\equiv r\ ({\rm mod}\
a)\big\}={\mathcal A_1}\bigcup {\mathcal A_2},
\end{equation}
where
$${\mathcal A_1}=\bigcup_{i=0}^{H-1}\Big\{{\rm prime}\ q\equiv r({\rm mod}\ a):
\frac{b+amn}{r'+ai}<q\le \frac{b+aln}{r'+ai}\Big\}$$ and
$${\mathcal A_2}=\Big\{{\rm
prime}\ q\equiv r({\rm mod}\ a): q\le \frac{b+aln}{r'+aH}\Big\}.$$

Since $H=\lfloor \frac{al-(l-m)r'}{a(l-m)}\rfloor $, we have
$al-(l-m)r'\ge aH(l-m)$. Thus for any positive integer  $i\le H$, we
obtain
\begin{align*}
\frac{b+amn}{r'+a(i-1)}-\frac{b+aln}{r'+ai}&=
\frac{a\big(b+aln-(l-m)r'n-ai(l-m)n\big)}{(r'+a(i-1))(r'+ai)}\\
&> \frac{a(aH(l-m)n-ai(l-m)n)}{(r'+a(i-1))(r'+ai)}\\
&= \frac{a^2(H-i)(l-m)n}{(r'+a(i-1))(r'+ai)}\ge 0.
\end{align*}
That is, for any positive integer  $i\le H$, one has
$$
\frac{b+amn}{r'+a(i-1)}> \frac{b+aln}{r'+ai}.
$$
It follows that for any integer $1\le i\le H$, the intersection of
$$\Big\{{\rm prime}\ q\equiv r({\rm mod}\ a):
\frac{b+amn}{r'+ai}<q\le \frac{b+aln}{r'+ai}\Big\}$$ and
$$\Big\{{\rm prime}\ q\equiv r({\rm mod}\ a): \frac{b+amn}{r'+a(i-1)}<q\le
\frac{b+aln}{r'+a(i-1)}\Big\}$$ is empty and ${\mathcal A_1}\cap
{\mathcal A_2}$ is empty too.

First using (\ref{eq:2.12}), and then using (\ref{eq:1.1}) and
(\ref{eq:1.2}), we obtain that if $l\le \frac{(a+r')m}{r'}$, then
\begin{align}\label{eq:2.13}
&\sum_{p\in P_{m,l}(n)\atop p\equiv r\pmod a}\log
p\\
\nonumber&=\sum_{p\in {\mathcal A_1}}\log p+\sum_{p\in {\mathcal
A_2}}\log p =\sum_{i=0}^{H-1}\sum_{p\equiv r\pmod a\atop
\frac{b+amn}{r'+ai}<p\le \frac{b+aln}{r'+ai}}\log p+\sum_{p\equiv
r\pmod a\atop
p\le \frac{b+aln}{r'+aH}}\log p\\
\nonumber&=\sum_{i=0}^{H-1}\bigg( \vartheta\Big(
\frac{b+aln}{r'+ai}; a,r\Big) -\vartheta\Big( \frac{b+amn}{r'+ai};
a,
r\Big)\bigg)+\vartheta\Big( \frac{b+aln}{r'+aH}; a, r\Big)\\
\nonumber&=\sum_{i=0}^{H-1}\frac{1}{\varphi(a)}
\frac{a(l-m)n}{r'+ai}+ \frac{1}{\varphi(a)}
\frac{aln}{r'+aH}+O\big(n\exp(-c_1\sqrt{\log
n})\big)\\
\nonumber&=\frac{an}{\varphi(a)}\Big( \sum_{i=0}^{H-1}
\frac{l-m}{r'+ai}+
 \frac{l}{r'+aH}\Big)+O\big(n\exp(-c_1\sqrt{\log
n})\big).
\end{align}

Evidently, one has $H=1$ and $\frac{l}{r'+a}=\frac{m}{r'}$ if
$l=\frac{(a+r')m}{r'}$. Hence
\begin{equation}\label{eq:2.14}
\frac{a}{\varphi(a)}\Big( \sum_{i=0}^{H-1}
\frac{l-m}{r'+ai}+\frac{l}{r'+aH}\Big)=\frac{al}{r'\varphi(a)}
\end{equation}
if $l=\frac{(a+r')m}{r'}$.

Now by (\ref{eq:2.3}), (\ref{eq:2.10}), (\ref{eq:2.13}) and
(\ref{eq:2.14}), we can deduce from $H=K_{r'}$ that
\begin{align}\label{eq:2.15}
\log L_{m, l}(n)&=\sum_{r\in R(a)}\sum_{p\in P_{m, l}(n)\atop
p\equiv r\pmod a}\log p
+O(\sqrt{n})\\
\nonumber&=\frac{an}{\varphi(a)}\Big(\sum_{r\in
R(a)}A_{r'}\Big)+O\big(n\exp(-c_1\sqrt{\log n})\big),
\end{align}
where $A_{r'}$ is defined as in (\ref{eq:1.3}).

Since $r'$ runs over $R(a)$ as $r$ does, we have
$$
\sum_{r\in R(a)} A_{r'}=\sum_{r=1\atop \gcd(a, r)=1}^a A_r.
$$
It then follows immediately from (2.15) that
$$\log L_{m, l}(n)=\frac{an}{\varphi(a)}\Big(\sum_{r=1\atop \gcd(a, r)=1}^a A_r\Big)+o(n)$$
as desired. This completes the proof of Theorem \ref{mainresult}.
\hfill$\Box$


\subsection*{Acknowledgment}
The authors would like to thank the anonymous referee for helpful
comments and suggestions.

\end{document}